\documentclass[11pt,onecolumn]{IEEEtran}
\usepackage{graphicx}
\usepackage{amsmath}
\usepackage{mathrsfs}
\usepackage{mathtools}
\usepackage{amssymb}
\usepackage{float}
\usepackage{url}
\usepackage{verbatim}
\usepackage{enumerate}
\usepackage{layout}
\usepackage{bm}

\usepackage[bottom=0.79in]{geometry}

\newtheorem{theorem}{Theorem}

\newtheorem{proposition}{Proposition}

\newtheorem{remark}{Remark}
\newcommand{\prob}{\ensuremath{\mathbb{P}}}
\newcommand{\naturals}{\ensuremath{\mathbb{N}}}
\newcommand{\Reals}{\ensuremath{\mathbb{R}}}

\newcommand{\set}{\ensuremath{\mathcal}}

\newcommand{\OneTo}[1]{[1:#1]}
\newcommand{\eqdef}{\triangleq} 
\newcommand{\card}[1]{|#1|}
\newcommand{\bigcard}[1]{\bigl|#1\bigr|}

\RequirePackage{bbm}

\newcommand{\pmfsub}[1]{{\smash{#1}\vphantom{XYZ}}}
\newcommand{\pmf}[1]{\mathsf{#1}}
\newcommand{\pmfOf}[1]{\pmf{P}_{\pmfsub{#1}}}  
\newcommand{\CondpmfOf}[2]{\pmf{P}_{\pmfsub{#1} | \pmfsub{#2}}}

\DeclareMathOperator{\Vertex}{\mathsf{V}}
\DeclareMathOperator{\Edge}{\mathsf{E}}
\DeclareMathOperator{\Independentset}{\set{I}}
\DeclareMathOperator{\Neighbors}{\set{N}}
\DeclareMathOperator{\Degree}{\text{d}}
\newcommand{\V}[1]{\Vertex(#1)}
\newcommand{\E}[1]{\Edge(#1)}
\newcommand{\indset}[1]{\Independentset(#1)}  
\newcommand{\Ngb}[1]{\Neighbors(#1)} 
\newcommand{\dgr}[1]{\Degree_{#1}}   

\DeclareMathOperator{\Entr}{H}
\DeclareMathOperator{\Info}{I}
\newcommand{\Hb}{\mathsf{H}_{\textnormal{b}}}

\newcommand{\Ent}[1]{\Entr(#1)}

\newcommand{\EntBin}[1]{\Hb(#1)}

\newcommand{\biggEntBin}[1]{\Hb\biggl(#1\biggr)}

\newcommand{\EntCond}[2]{\Entr(#1 | \kern0.1em #2)}
\newcommand{\bigEntCond}[2]{\Entr\bigl(#1 | \kern0.1em #2\bigr)}
\newcommand{\BigEntCond}[2]{\Entr\Bigl(#1 \kern-0.1em \bigm| \kern-0.1em #2 \Bigr)}
\newcommand{\biggEntCond}[2]{\Entr\biggl(#1 \kern-0.1em \Bigm| \kern-0.1em #2 \biggr)}
\newcommand{\BiggEntCond}[2]{\Entr\Biggl(#1 \kern-0.1em \biggm| \kern-0.1em #2 \Biggr)}

\newcommand{\biggMInfo}[2]{\Info\biggl(#1; \kern0.1em #2 \biggr)}
\newcommand{\BiggMInfo}[2]{\Info\Biggl(#1; \kern0.1em #2 \Biggr)}

\newcommand{\KLDivBin}[2]{\mathsf{D}_{\textnormal{b}}(#1 \kern0.1em \| \kern0.1em #2)} 
\newcommand{\bigKLDivBin}[2]{\mathsf{D}_{\textnormal{b}}\bigl(#1 \kern0.1em \| \kern0.1em #2 \bigr)}
\newcommand{\BigKLDivBin}[2]{\mathsf{D}_{\textnormal{b}}\Bigl(#1 \kern0.1em \| \kern0.1em #2 \Bigr)}
\newcommand{\biggKLDivBin}[2]{\mathsf{D}_{\textnormal{b}}\biggl(#1 \kern0.1em \| \kern0.1em #2 \biggr)}
\newcommand{\BiggKLDivBin}[2]{\mathsf{D}_{\textnormal{b}}\Biggl(#1 \kern0.1em \| \kern0.1em #2 \Biggr)}
\newcommand{\bigKLDiv}[2]{\mathsf{D}\bigl(#1 \kern0.1em \| \kern0.1em #2 \bigr)}
\newcommand{\BigKLDiv}[2]{\mathsf{D}\Bigl(#1 \kern-0.2em \bigm\| \kern-0.2em #2 \Bigr)}
\newcommand{\biggKLDiv}[2]{\mathsf{D}\biggl(#1 \kern-0.2em \Bigm\| \kern-0.2em #2 \biggr)}
\newcommand{\BiggKLDiv}[2]{\mathsf{D}\Biggl(#1 \kern-0.1em \biggm\| \kern-0.1em #2 \Biggr)}

\newcommand{\I}[1]{{\mathbbm 1}\{#1\}}

\begin{document}
\thispagestyle{empty}
\setcounter{page}{1}
\setlength{\baselineskip}{1.15\baselineskip}

\title{\huge{A Generalized Information-Theoretic Approach for Bounding the
Number of Independent Sets in Bipartite Graphs}\\[0.2cm]}
\author{Igal Sason
\thanks{
I. Sason is with the Andrew \& Erna Faculty of Electrical Engineering,
Technion - Israel Institute of Technology, Haifa 3200003, Israel
(e-mail: sason@ee.technion.ac.il).}}

\maketitle
\thispagestyle{empty}

\vspace*{-0.4cm}
\begin{abstract}
This paper studies the problem of upper bounding the number of independent sets in a graph,
expressed in terms of its degree distribution. For bipartite regular graphs, Kahn (2001) established
a tight upper bound using an information-theoretic approach, and he also conjectured an upper bound
for general graphs. His conjectured bound was recently proved by Sah et al. (2019), using different
techniques not involving information theory. The main contribution of this work is the extension of
Kahn’s information-theoretic proof technique to handle irregular bipartite graphs. In particular, when
the bipartite graph is regular on one side, but it may be irregular in the other, the extended entropy-based
proof technique yields the same bound that was conjectured by Kahn (2001) and proved by Sah et al. (2019).
\end{abstract}

{\bf{Keywords}}: {\small Shannon entropy, Shearer's lemma, counting, independent sets, graphs.}

\vspace*{-0.2cm}
\section{Introduction}
\label{section: introduction}

The Shannon entropy and other classical information measures serve as a powerful tool in various
combinatorial and graph-theoretic applications (e.g., \cite{BabuR14}--\cite{Friedgut04},
\cite{Galvin14}--\cite{Radhakrishnan01} and \cite{Spencer85})
such as the method of types, applications of Shearer's lemma,
sub and supermodularity properties of information measures and their applications,
entropy-based proofs of Moore bound for irregular graphs, Bregman's theorem on the permanent
of square matrices with binary entries, Spencer's theorem in discrepancy theory, etc.

The enumeration of discrete structures that satisfy certain local constraints, and
particularly the enumeration of independent sets in graphs is of interest in discrete mathematics.
Many important structures can be modeled by independent sets in a graph, i.e., subsets of vertices
in a graph where none of them are connected by an edge. For example, if a graph models some kind
of incompatibility, then an independent set in this graph represents a mutually compatible collection.
Upper bounding the number of independent sets in a regular graph was motivated in \cite{Alon91} by a
conjecture which has several applications in combinatorial group theory. A survey paper on upper bounding
the number of independent sets in graphs, along with some of their applications, is provided in \cite{Samotij15}.
The problem of counting independent sets in graphs received, in general, significant attention in the
literature of discrete mathematics over the last three decades, and also in the information theory literature
(\cite{MadimanT_ISIT07, MadimanT_IT10}).

A tight upper bound on the number of independent sets in finite
and undirected general graphs was proved in the special setting of bipartite regular graphs
in \cite{Kahn01}, and it was conjectured there to hold for general (irregular) graphs (2001,
see Conjecture~4.2 in \cite{Kahn01}).
A decade later (2010), it was extended in \cite{Zhao10} to regular graphs (that are not
necessarily bipartite); a year later (2011), it was proved in \cite{GalvinZ11} for graphs with
small degrees (up to~5). Finally, this conjecture was recently (2019) proved in general \cite{SahSaStZhao19},
by utilizing a new approach. The reader is referred to \cite{SahSaStZhao19b} for an announcement
on the solution of this conjecture as a frustrating combinatorial problem for two decades, along with the history
and ramifications of this problem, and some reflections of the authors on their work in \cite{SahSaStZhao19}.

The recently-introduced proof of the conjecture for {\em general} undirected graphs \cite{SahSaStZhao19}
uses an induction on the number of vertices in a graph, and it obtains a recurrence inequality whose
derivation involves some judicious applications of H\"{o}lder's inequality (see
Sections~2 and 4 in \cite{SahSaStZhao19}). It settled for the first time
Conjecture~4.2 in \cite{Kahn01} by an interesting approach, which is unrelated to information theory.
The possibility of generalizing the information-theoretic
proof in \cite{Kahn01} to irregular bipartite graphs was left in \cite{SahSaStZhao19} as an open issue.
It should be noted that by proving Kahn's conjecture for irregular bipartite graphs, this readily enables to
extend the proof to general undirected graphs (by invoking Zhao's inequality, see Lemma~3 in \cite{GalvinZ11}).

The main contribution of this work is the extension of Kahn’s information-theoretic proof
technique for bipartite regular graphs \cite{Kahn01} to handle irregular bipartite graphs. In particular, when
the bipartite graph is regular on one side, but it may be irregular in the other, the extended entropy-based
proof technique yields the same bound that was conjectured by Kahn \cite{Kahn01} and proved by Sah et al. \cite{Zhao10}.

The structure of the paper is as follows: Section~\ref{section: preliminaries} provides
preliminaries and notation that are essential for the analysis in this paper.
Section~\ref{section: scientific merit}
explains (in more details) the scientific merit and contributions of the present work; for the sake
of causal presentation, we provide these explanations after Section~\ref{section: preliminaries}.
Sections~\ref{section: An IT proof of Theorem 3} and~\ref{section: An IT proof of Zhao's inequality}
are the core of this work.

\section{Preliminaries and Notation}
\label{section: preliminaries}

We provide in this section the notation and preliminary material which is essential for the presentation in
this paper.

\subsection{Notation and Basic Properties of the Entropy}
\label{subsection: notation}

The following notation is used in the present paper:
\begin{itemize}
\item $\naturals \eqdef \{1, 2, \ldots \}$ denotes the set of natural numbers.
\item $X^n \eqdef (X_1, \ldots X_n)$ denotes an $n$-dimensional random vector
of discrete random variables, having a joint probability mass function (PMF)
that is denoted by $\pmfOf{X^n}$.
\item For every $n \in \naturals$, $\OneTo{n} \eqdef \{1, \ldots, n\}$;
\item $X_{\set{S}} \eqdef (X_i)_{i \in \set{S}}$ is a random vector for
an arbitrary nonempty subset $\set{S} \subseteq \OneTo{n}$; if
$\set{S} = \emptyset$, then conditioning on $X_{\set{S}}$ is void.
\item $\I{E}$ denotes the indicator of an event $E$; i.e., it is equal to~1
if this event is satisfied, and it is zero otherwise.
\item Let $X$ be a discrete random variable that takes its values on a set $\set{X}$,
and let $\pmfOf{X}$ be the PMF of $X$. The {\em Shannon entropy} of $X$ is given by
\begin{IEEEeqnarray}{rCl}
\label{eq: entropy}
\Ent{X} \eqdef -\sum_{x \in \set{X}} \pmfOf{X}(x) \, \log \pmfOf{X}(x),
\end{IEEEeqnarray}
where throughout this paper, we take all logarithms to base~2.
\item For $p \in [0,1]$,
\begin{IEEEeqnarray}{rCl}
\label{eq2: EntBin}
\EntBin{p} \eqdef -p \log p -(1-p) \log(1-p),
\end{IEEEeqnarray}
where $\EntBin{\cdot}$ is the {\em binary entropy function}.
By continuous extension, the convention $0 \log 0 = 0$ is used.
\item Let $X$ and $Y$ be discrete random variables with a joint PMF $\pmfOf{XY}$, and having a marginal
PMF of $X$ given $Y$ which is denoted by $\CondpmfOf{X}{Y}$. Let $X$ and $Y$ take their values in the
sets $\set{X}$ and $\set{Y}$, respectively. The {\em conditional entropy} of $X$ given $Y$
is defined as
\begin{IEEEeqnarray}{rCl}
\label{eq: conditional entropy 1}
\EntCond{X}{Y} & \eqdef & -\sum_{(x,y) \in \set{X} \times \set{Y}} \pmfOf{XY}(x,y) \log \CondpmfOf{X}{Y}(x|y) \\[0.1cm]
\label{eq: conditional entropy 1.5}
&=& \sum_{y \in \set{Y}} \pmfOf{Y}(y) \, \EntCond{X}{Y=y},
\end{IEEEeqnarray}
and
\begin{IEEEeqnarray}{rCl}
\EntCond{X}{Y} &=& \Ent{X,Y} - \Ent{Y}.  \label{eq: conditional entropy 2}
\end{IEEEeqnarray}
\end{itemize}

This paper relies on the following basic properties of the Shannon entropy:
\begin{itemize}
\item Entropies and conditional entropies of discrete random vectors are nonnegative.
\item If $\set{X}$ is a finite set, then
\begin{IEEEeqnarray}{rCl}
\label{eq: UB Ent}
\Ent{X} \leq \log \card{\set{X}},
\end{IEEEeqnarray}
with equality in \eqref{eq: UB Ent} if and only if $X$ is equiprobable over the set $\set{X}$.
\item Conditioning cannot increase the entropy, i.e.,
\begin{IEEEeqnarray}{rCl}
\label{eq: conditioning reduces entropy}
\EntCond{X}{Y} \leq \Ent{X},
\end{IEEEeqnarray}
with equality in \eqref{eq: conditioning reduces entropy} if and only if $X$ and $Y$
are independent.
\item Generalizing \eqref{eq: conditional entropy 2}
to $n$-dimensional random vectors gives the {\em chain rule for the entropy}:
\begin{IEEEeqnarray}{rCl}
\label{eq1: chain rule}
\Ent{X^n} &=& \Ent{X_1} + \EntCond{X_2}{X_1} + \ldots + \EntCond{X_n}{X_1, \ldots, X_{n-1}} \\[0.1cm]
\label{eq2: chain rule}
&=& \sum_{i=1}^n \EntCond{X_i}{X^{i-1}}.
\end{IEEEeqnarray}
\item The following {\em subadditivity property of the entropy} is implied by
\eqref{eq: conditioning reduces entropy} and \eqref{eq2: chain rule}:
\begin{IEEEeqnarray}{rCl}
\label{eq: subadditivity Ent}
\Ent{X^n} \leq \sum_{i=1}^n \Ent{X_i},
\end{IEEEeqnarray}
with equality in \eqref{eq: subadditivity Ent} if and only if $X_1, \ldots, X_n$ are
independent random variables.
\end{itemize}

\subsection{Shearer's Lemma}
\label{subsection: Shearer's Lemma}
Shearer's lemma extends the subadditivity property \eqref{eq: subadditivity Ent} of the entropy.
Due to its simplicity and usefulness in this paper (and elsewhere), we state and
prove it here.
\begin{proposition}[Shearer's Lemma, \em{\cite{ChungGFS86}}]
\label{proposition: Shearer's Lemma}
Let $X_1, \ldots, X_n$ be discrete random variables, and let the sets
$\set{S}_1, \ldots, \set{S}_m \subseteq \OneTo{n}$ include every
element $i \in \OneTo{n}$ in {\em at least} $k \geq 1$ of these
subsets. Then,
\begin{IEEEeqnarray}{rCl}
\label{eq: Shearer's Lemma}
k \Ent{X^n} \leq \sum_{j=1}^m \Ent{X_{\set{S}_j}}.
\end{IEEEeqnarray}
\end{proposition}
As a special case of \eqref{eq: Shearer's Lemma}, setting $\set{S}_i \eqdef \{i\}$
as singletons for all $i \in \OneTo{n}$ gives \eqref{eq: subadditivity Ent}
by having $k=1$ and $m=n$.

\begin{IEEEproof}
Let $\set{S} = \{i_1, \ldots, i_\ell\}$ with $1 \leq i_1 < \ldots < i_\ell \leq n$.
By invoking the chain rule in this order,
\begin{IEEEeqnarray}{rCl}
\Ent{X_{\set{S}}} &=& \Ent{X_{i_1}} + \EntCond{X_{i_2}}{X_{i_1}} + \ldots +
\EntCond{X_{i_\ell}}{X_{i_1}, \ldots, X_{i_{\ell-1}}} \nonumber \\[0.1cm]
&\geq& \sum_{i \in \set{S}} \EntCond{X_i}{X^{i-1}} \nonumber \\
&=& \sum_{i=1}^n \I{i \in \set{S}} \, \EntCond{X_i}{X^{i-1}},
\label{eq: conditioning reduces entropy 2}
\end{IEEEeqnarray}
where the last inequality holds since an additional conditioning cannot increase
the entropy (i.e., $\EntCond{X}{Y,Z} \leq \EntCond{X}{Y}$ for all $X,Y$ and $Z$).
By assumption, for $i \in \OneTo{n}$,
\begin{IEEEeqnarray}{rCl}
\label{eq: degree}
\sum_{j=1}^m  \I{i \in \set{S}_j} \geq k,
\end{IEEEeqnarray}
since the number of subsets $\{\set{S}_j\}_{j=1}^m$ that include $i$ as an element is at least $k$. Consequently,
it follows that
\begin{IEEEeqnarray}{rCl}
\sum_{j=1}^m \Ent{X_{\set{S}_j}} & \geq &
\sum_{j=1}^m \sum_{i=1}^n \bigl\{ \I{i \in \set{S}_j} \, \EntCond{X_i}{X^{i-1}} \bigr\} \label{151220a1} \\
&=& \sum_{i=1}^n \Biggl\{ \sum_{j=1}^m  \I{i \in \set{S}_j} \, \EntCond{X_i}{X^{i-1}} \Biggr\} \label{151220a2} \\
&\geq& k \sum_{i=1}^n  \EntCond{X_i}{X^{i-1}} \label{151220a3} \\[0.1cm]
&=& k \Ent{X^n},  \label{151220a4}
\end{IEEEeqnarray}
where \eqref{151220a1} follows from \eqref{eq: conditioning reduces entropy 2};
\eqref{151220a2} holds by interchanging the order of summation;
\eqref{151220a3} holds by \eqref{eq: degree};
and \eqref{151220a4} holds by the chain rule for the Shannon entropy.
\end{IEEEproof}

\begin{remark}
\label{remark: Sheaerer's lemma}
Inequality \eqref{eq: Shearer's Lemma} holds even if the sets $\set{S}_1, \ldots, \set{S}_m$ are
not necessarily included in $\OneTo{n}$. To verify it, define the subsets
$\set{S}'_j \eqdef \set{S}_j \cap [1:n]$ for all $j \in \OneTo{m}$. The subsets
$\set{S}'_1, \ldots, \set{S}'_m$ are all included in $\OneTo{n}$, and every element
$i \in \OneTo{n}$ continues to be included in at least $k \geq 1$ of these subsets. Hence,
Proposition~\ref{proposition: Shearer's Lemma} can be applied to the subsets
$\set{S}'_1, \ldots, \set{S}'_m$. By the monotonicity property of the entropy, the inclusion
$\set{S}'_j \subseteq \set{S}_j$ implies that $\Ent{X_{\set{S}'_j}} \leq \Ent{X_{\set{S}_j}}$
for all $j \in \OneTo{m}$, which then yields the statisfiability of \eqref{eq: Shearer's Lemma}.
\end{remark}

\begin{remark}
\label{remark2: Shearer's lemma}
A generalized inequality which extends both Shearer's lemma and Han's inequality is provided
in \cite[Proposition~1]{MadimanT_IT10}.
\end{remark}

\vspace*{0.1cm}
Shearer's lemma and some of its variants (see \cite{Friedgut04} and
\cite{GavinkyLSS14}) have been successfully applied in various occasions
(see, e.g., \cite{Friedgut04, Galvin14, GavinkyLSS14, Kahn01, Kahn02, Radhakrishnan01}).
Shearer's lemma is also instrumental in this paper.

\subsection{Graphs, Independent Sets, and Tensor Products}
\label{subsection: Graphs}

Let $G$ be an undirected graph, and let $\V{G}$ and $\E{G}$ denote respectively the sets
of vertices and edges in $G$.

A graph $G$ is called {\em $d$-regular} if the degree of all vertices in $\V{G}$ is equal
to $d$. Otherwise, if the graph $G$ is not $d$-regular for some $d \in \naturals$, then
$G$ is an {\em irregular} graph.

A graph is called {\em bipartite} if it has two types of vertices, and an edge cannot
connect vertices of the same type; we refer to the vertices of a bipartite graph $G$
as left and right vertices.

A graph $G$ is called {\em complete} if every vertex $v \in \V{G}$ is connected to all the
other vertices in $\V{G} \setminus \{v\}$ (and not to itself); similarly, a bipartite graph
is called complete if every vertex is connected to all the vertices of the other type in
the graph. A complete $(d-1)$-regular graph is denoted by $K_d$, having a number of vertices
$\bigcard{\V{K_d}} = d$, and a number of edges $\bigcard{\E{K_d}} = \tfrac12 \, d(d-1)$.
Likewise, a complete $d$-regular bipartite graph is denoted by $K_{d,d}$, having a number
of vertices $\bigcard{\V{K_{d,d}}} = 2d$ (i.e., $d$ vertices of each of the two types),
and a number of edges $\bigcard{\E{K_{d,d}}} = d^2$.

An {\em independent set} of an undirected graph $G$ is a subset of its vertices such that
none of the vertices in this subset are adjacent (i.e., none of them are joined by an edge).
Let $\indset{G}$ denote the set of all the independent sets in $G$, and let $\bigcard{\indset{G}}$
denote the number of independent sets in $G$. Similarly to
\cite{Alon91, GalvinZ11, Galvin14, Kahn01, MadimanT_ISIT07, MadimanT_IT10, Radhakrishnan01, SahSaStZhao19, Samotij15, Zhao10}
(and references therein), our work considers the question of how many independent sets
can $G$ have.

The {\em tensor product} $G \times H$ of two graphs $G$ and $H$ is a graph
such that the following holds:
\begin{itemize}
\item The vertex set of $G \times H$ is the Cartesian product $\V{G} \times \V{H}$;
\item Two vertices $(g,h), (g', h') \in \V{G \times H}$ are adjacent if and only if
$g$ is adjacent to $g'$, and $h$ is adjacent to $h'$, i.e., $(g,g') \in \E{G}$ and
$(h,h') \in \E{H}$. This is denoted by $\bigl( (g,h), (g', h') \bigr) \in \E{G \times H}$.
\end{itemize}
In general, the following identities hold:
\begin{IEEEeqnarray}{rCl}
& \bigcard{\V{G \times H}} = \bigcard{\V{G}} \, \bigcard{\V{H}}, \label{eq: vertices of tensor product} \\[0.1cm]
& \bigcard{\E{G \times H}} = 2 \, \bigcard{\E{G}} \, \bigcard{\E{H}}.  \label{eq: edges of tensor product}
\end{IEEEeqnarray}

By the definition of a complete $d$-regular graph $K_d$, the graph $K_2$
is specialized to two vertices that are connected by an edge. Let us label the
two vertices in $K_2$ by 0 and~1. For a graph $G$, the tensor product
$G \times K_2$ is a bipartite graph, called the {\em bipartite double cover}
of $G$, where the set of vertices in $G \times K_2$ is given by
\begin{IEEEeqnarray}{rcl}
\label{eq: vertices of bipartite double cover}
\V{G \times K_2} = \bigl\{ (v,i) : v \in \V{G}, \, i \in \{0,1\} \bigr\},
\end{IEEEeqnarray}
and its set of edges is given by
\begin{IEEEeqnarray}{rCl}
\label{eq: edges of bipartite double cover}
\E{G \times K_2} = \bigl\{ \bigl( (u,0), (v,1) \bigr) : (u,v) \in \E{G} \bigr\}.
\end{IEEEeqnarray}
Hence, every edge $e = (u,v) \in \E{G}$ is mapped into the two edges
$\bigl( (u,0), (v,1) \bigr) \in \E{G \times K_2}$ and
$\bigl( (v,0), (u,1) \bigr) \in \E{G \times K_2}$
(since the graph $G$ is undirected). This implies that the numbers of vertices
and edges in $G \times K_2$ are doubled in comparison to their respective numbers
in $G$; moreover, every edge in $G$, which connects a pair of
vertices of specified degrees, is mapped into two edges in $G \times K_2$ where
each of these two edges connects a pair of vertices of the same specified degrees.

\subsection{Upper Bounds on the Number of Independent Sets}
\label{subsection: UB - Independent sets}

The present subsection introduces the relevant results to this paper.
The following theorem provides a tight upper bound on the number of independent
sets in bipartite regular graphs, and its derivation in \cite{Kahn01} makes
a clever use of Shearer's lemma (Proposition~\ref{proposition: Shearer's Lemma}).
\begin{theorem}[Kahn 2001, \cite{Kahn01}]
\label{theorem: Kahn 2001}
If $G$ is a bipartite $d$-regular graph with $n$ vertices, then
\begin{IEEEeqnarray}{rCl}
\label{eq: Kahn 2001}
\bigcard{\indset{G}} \leq \bigl( 2^{d+1} - 1 \bigr)^{\frac{n}{2d}}.
\end{IEEEeqnarray}
Furthermore, if $n$ is an even multiple of $d$, then the upper bound in
the right side of \eqref{eq: Kahn 2001} is tight, and it is obtained
by a disjoint union of $\frac{n}{2d}$ complete $d$-regular bipartite
graphs $(K_{d,d})$.
\end{theorem}

Kahn's result was later extended by Zhao \cite{Zhao10} for a general $d$-regular
graph via a brilliant combinatorial reduction to the setting of $d$-regular
bipartite graphs, which proved \cite[Conjecture~4.1]{Kahn01}.
\begin{theorem}[Zhao 2010, \cite{Zhao10}]
\label{theorem: Zhao10}
The upper bound on the number of independent sets in \eqref{eq: Kahn 2001}
continues to hold for all $d$-regular graphs with $n$ vertices.
\end{theorem}

Recently, Sah {\em et al.} \cite{SahSaStZhao19} proved Kahn's conjecture in
\cite[Conjecture~4.2]{Kahn01} (made eighteen years earlier) for an upper bound on
the number of independent sets in a general undirected graph with no isolated
vertices. The proof in \cite{SahSaStZhao19} is combinatorial, and it extends
the result in Theorem~\ref{theorem: Zhao10} as follows.
\begin{theorem}[Sah et al. 2019, \cite{SahSaStZhao19}]
\label{theorem: Sah et al., 2019}
Let $G$ be an undirected graph without isolated vertices or multiple edges
connecting any pair of vertices. Let $\dgr{v}$ denote the degree of a vertex
$v \in \V{G}$. Then,
\begin{IEEEeqnarray}{rCl}
\label{eq: Sah et al.}
\bigcard{\indset{G}} \leq \prod_{(u,v) \in \E{G}}
(2^{\dgr{u}} + 2^{\dgr{v}} - 1)^{\frac1{\dgr{u} \dgr{v}}}
\end{IEEEeqnarray}
with an equality if $G$ is a disjoint union of complete bipartite graphs.
\end{theorem}

Let $K_{\dgr{u}, \dgr{v}}$ be a complete bipartite graph where the degrees
of its left and right vertices are equal to $\dgr{u}$ and $\dgr{v}$, respectively.
Then, the number of independent sets in such a complete bipartite graph is
equal to
\begin{IEEEeqnarray}{rcl}
\label{eq: complete bipartite graph}
\bigcard{\indset{K_{\dgr{u}, \dgr{v}}}} = 2^{\dgr{u}} + 2^{\dgr{v}} - 1
\end{IEEEeqnarray}
since every subset of the left vertices, as well as every subset of the right
vertices, forms an independent set of $K_{\dgr{u}, \dgr{v}}$; on the other hand, any subset
which contains both left and right vertices is not an independent set (since
the bipartite graph $K_{\dgr{u}, \dgr{v}}$ is complete). Note that the
substraction by~1 in the right side of \eqref{eq: complete bipartite graph} is because,
in the counting of the number of subsets of left vertices ($2^{\dgr{v}}$)
or right vertices ($2^{\dgr{u}}$), the empty set is counted twice. Hence,
\eqref{eq: Sah et al.} can be rewritten in an equivalent form as
\begin{IEEEeqnarray}{rCl}
\label{eq2: Sah et al.}
\bigcard{\indset{G}} \leq \prod_{(u,v) \in \E{G}}
\bigcard{\indset{K_{\dgr{u}, \dgr{v}}}}^{\frac1{\dgr{u} \dgr{v}}}.
\end{IEEEeqnarray}
Since $\E{K_{\dgr{u}, \dgr{v}}} = \dgr{u} \dgr{v}$, it follows that
the bound in \eqref{eq: Sah et al.} (or \eqref{eq2: Sah et al.}) is
achieved by the complete bipartite graph $K_{\dgr{u}, \dgr{v}}$.
More generally, the bound is achieved by a disjoint finite union
of such complete bipartite graphs since the number of independent sets
in a disjoint union of graphs is equal to the product of the number
of independent sets in each of these component graphs.

For the extension of the validity of Theorem~\ref{theorem: Kahn 2001} to
Theorem~\ref{theorem: Zhao10}, obtained by relaxing the requirement
that the graph is bipartite, the following inequality was introduced by
Zhao for every finite graph $G$ \cite[Lemma~2.1]{Zhao10}:
\begin{IEEEeqnarray}{rCl}
\label{eq: Zhao's inequality}
\bigcard{\indset{G}}^2 \leq \bigcard{\indset{G \times K_2}},
\end{IEEEeqnarray}
which relates the number of independent sets in a graph to the number
of independent sets in the bipartite double cover of this graph.

The transition from Theorem~\ref{theorem: Kahn 2001}
to Theorem~\ref{theorem: Zhao10}, as introduced in \cite{Zhao10},
is a one-line proof. Let $G$ be a $d$-regular graph
with $n$ vertices, then $G \times K_2$ is $d$-regular bipartite
graph with $2n$ vertices. Hence, \eqref{eq: Kahn 2001} and
\eqref{eq: Zhao's inequality} give that
\begin{IEEEeqnarray}{rCl}
\label{eq: Zhao's proof}
\bigcard{\indset{G}}^2 \leq \bigcard{\indset{G \times K_2}}
\leq (2^{d+1}-1)^{\frac{2n}{2d}},
\end{IEEEeqnarray}
and taking the square root of the leftmost and rightmost sides of \eqref{eq: Zhao's proof}
implies that \eqref{eq: Kahn 2001} continues to hold even when the
regular graph $G$ is not necessarily bipartite.

\section{Scientific Merit and Contributions of the Present Work}
\label{section: scientific merit}

After introducing Shearer's lemma (Proposition~\ref{proposition: Shearer's Lemma}) and
Theorems~\ref{theorem: Kahn 2001}--\ref{theorem: Sah et al., 2019}, we address the scientific
merit and contributions of the present work in more details (in comparison to the introduction
in Section~\ref{section: introduction}).

Theorem~\ref{theorem: Sah et al., 2019} was recently proved in \cite{SahSaStZhao19} (see also \cite{SahSaStZhao19b})
for general graphs, without relying on information theory.
The motivation of our work is rooted in the following sentences from \cite[p.~174]{SahSaStZhao19}:
\par
{\em Kahn's proof \cite{Kahn01} of the bipartite case of Theorem~\ref{theorem: Kahn 2001} made
clever use of Shearer's entropy inequality \cite{ChungGFS86}.
It remains unclear how to apply Shearer's inequality in a lossless way in the irregular case, despite
previous attempts to do so, e.g., \cite[Section~3]{MadimanT_ISIT07} and \cite[Section~5.C]{MadimanT_IT10}}.

The present paper gives an information-theoretic proof of Theorem~\ref{theorem: Sah et al., 2019}
in a setting where the bipartite graph is regular on one side (i.e., the vertices on the
other side of the bipartite graph can be irregular, and have arbitrary degrees). Its contributions
are as follows:
\begin{itemize}
\item Section~\ref{section: An IT proof of Theorem 3} provides
a (non-trivial) extension of the proof of \eqref{eq: Kahn 2001},
from regular bipartite graphs \cite{Kahn01} to general bipartite graphs.
It leads to an upper bound on the number of independent sets which is
in general looser than the bound in \eqref{eq: Sah et al.} (or its
equivalent form in \eqref{eq2: Sah et al.}).
However, for the family of bipartite graphs that are regular on one side of
the graph, our bound in \eqref{eq1: UB cardinality} coincides with the bound
in \eqref{eq: Sah et al.}.
The main deviation from Kahn's information-theoretic proof in \cite{Kahn01}
is that we allow here the bipartite graph to be irregular. This generalization
is not trivial in the sense that it requires a more careful analysis and a slightly
more complicated version of Shearer’s Lemma (see Remark~\ref{remark: Sheaerer's lemma}).
It should be noted, however, that the suggested proof does follow the same recipe of
Kahn's proof in \cite{Kahn01} with some further complications that arise from the
non-regularity of the bipartite graphs.

\item A variant of the proof of Zhao's inequality \eqref{eq: Zhao's inequality}
(\cite[Section~2]{Zhao10}) is provided in
Section~\ref{section: An IT proof of Zhao's inequality}.
\end{itemize}

It is interesting to note that the observation that \eqref{eq: Sah et al.} can
be extended from (undirected) bipartite graphs to general graphs, by utilizing
\eqref{eq: Zhao's inequality}, was made in \cite[Lemma~3]{GalvinZ11}.
However, a computer-assisted proof of \eqref{eq: Sah et al.} was restricted
there to graphs whose maximal degrees are at most~5 (see
\cite[Theorem~2]{GalvinZ11}).

\section{An Information-Theoretic Proof of Theorem~\ref{theorem: Sah et al., 2019} for a Family of Bipartite Graphs}
\label{section: An IT proof of Theorem 3}

The core of the proof of Theorem~\ref{theorem: Sah et al., 2019} is proving
\eqref{eq: Sah et al.} for an undirected bipartite graph. We provide here
an extension of the entropy-based proof by Kahn \cite{Kahn01} from bipartite $d$-regular
graphs to general bipartite graphs, and then we prove \eqref{eq: Sah et al.} for the family
of bipartite graphs that are regular on one side. As it is explained in Section~\ref{section: scientific merit},
the proof in the present section follows the same recipe of Kahn's proof in \cite{Kahn01} with
some complications that arise from the non-regularity of the bipartite graphs. The following proof
deviates from the proof in \cite{Kahn01} already at its starting point, by a proper adaptation
of the proof technique to the general setting of irregular bipartite graphs, followed by a bit
more complicated usage of Shearer's lemma and a more involved analysis.

Consider first a general bipartite graph $G$ with a number of vertices $\card{\V{G}} = n$, and where none
of its vertices is isolated. Label them by the elements of $\OneTo{n}$. Let $\set{L}$
and $\set{R}$ be the vertices of the two types in $\V{G}$ (called, respectively, the
left and right vertices in $G$), so $\V{G} = \set{L} \cup \set{R}$ is a disjoint union.
Let $\set{D}_{\mathrm{L}}$ and $\set{D}_{\mathrm{R}}$ be, respectively, the sets of all
possible degrees of vertices in $\set{L}$ and $\set{R}$. For all $d \in \set{D}_{\mathrm{L}}$,
let $\set{L}_d$ be the set of vertices in $\set{L}$ having degree $d$,
and let $\set{R}_d$ be the set of vertices in $\set{R}$ that are adjacent to vertices
in $\set{L}_d$ (note that the vertices in $\set{R}_d$ are not necessarily those vertices
in $\set{R}$ having degree $d$, so the definitions of $\set{L}_d$ and $\set{R}_d$ differ, i.e.,
they are not similar up to a replacement of left vertices of degree $d$ with right vertices
of the same degree). Then,
\begin{IEEEeqnarray}{rCl}
\label{eq: unions}
\set{L} = \bigcup_{d \in \set{D}_{\mathrm{L}}} \, \set{L}_d,
\quad \set{R} = \bigcup_{d \in \set{D}_{\mathrm{L}}} \, \set{R}_d,
\end{IEEEeqnarray}
where the first equality in \eqref{eq: unions} is (by definition) a union of pairwise disjoint sets.

Let $\set{S} \in \indset{G}$ be an independent set in $G$, which is selected uniformly at
random from $\indset{G}$, and let $X^n \eqdef (X_1, \ldots, X_n)$ be given by
\begin{IEEEeqnarray}{rCl}
\label{eq: indicators}
X_i \eqdef \I{i \in \set{S}}, \quad i \in \OneTo{n},
\end{IEEEeqnarray}
so the binary random vector $X^n$ indicates which of the $n$ vertices in $\V{G}$ belongs to the
randomly selected independent set $\set{S}$. Since $\set{S}$ is equiprobable in $\indset{G}$, we have
\begin{IEEEeqnarray}{rCl}
\label{eq: Ent01}
\Ent{X^n} = \log \, \bigcard{\indset{G}}.
\end{IEEEeqnarray}
Let $X_{\set{L}} = (X_i)_{i \in \set{L}}$ and $X_{\set{R}} = (X_i)_{i \in \set{R}}$ be used
as a shorthand. Then,
\begin{IEEEeqnarray}{rCl}
\label{eq: Ent02}
\Ent{X^n} &=& \Ent{X_{\set{L}}, X_{\set{R}}} \\[0.1cm]
\label{eq: Ent03}
&=& \Ent{X_{\set{L}}} + \EntCond{X_{\set{R}}}{X_{\set{L}}} \\[0.1cm]
\label{eq: Ent04}
& \leq & \sum_{d \in \set{D}_{\mathrm{L}}} \Ent{X_{\set{L}_d}} + \EntCond{X_{\set{R}}}{X_{\set{L}}} \\
\label{eq: Ent05}
& \leq & \sum_{d \in \set{D}_{\mathrm{L}}} \Ent{X_{\set{L}_d}}
+ \sum_{d \in \set{D}_{\mathrm{L}}} \EntCond{X_{\set{R}_d}}{X_{\set{L}}} \\
\label{eq: Ent06}
&=& \sum_{d \in \set{D}_{\mathrm{L}}} \bigl\{ \Ent{X_{\set{L}_d}}
+ \EntCond{X_{\set{R}_d}}{X_{\set{L}}} \bigr\},
\end{IEEEeqnarray}
where inequalities \eqref{eq: Ent04} and \eqref{eq: Ent05} hold by the subadditivity of the entropy,
and due to \eqref{eq: unions}. It should be noted that although the first summand in the right side
of \eqref{eq: Ent06} is an entropy of $X_{\set{L}_d}$, the conditioning on $X_{\set{L}}$ (rather
than just on $X_{\set{L}_d}$) in the second term leads to a stronger upper bound on $\Ent{X^n}$
(since $\set{L}_d \subseteq \set{L}$, and conditioning reduces the entropy). This is essential for
the continuation of the proof (see \eqref{eq: Ent08}).

We next upper bound the two summands in the right side of \eqref{eq: Ent06}, starting with
the conditional entropy. By invoking the subadditivity property of the entropy, for every
$d \in \set{D}_{\mathrm{L}}$,
\begin{IEEEeqnarray}{rCl}
\label{eq: Ent07}
\EntCond{X_{\set{R}_d}}{X_{\set{L}}} \leq \sum_{r \in \set{R}_d} \EntCond{X_r}{X_{\set{L}}}.
\end{IEEEeqnarray}
For every $r \in \set{R}_d$, let $\Ngb{r}$ be the set of all the vertices that are adjacent to the
vertex $r$. Since the graph $G$ is bipartite, we have $\Ngb{r} \subseteq \set{L}$ (but, in general,
$\Ngb{r} \not\subseteq \set{L}_d$), and consequently
\begin{IEEEeqnarray}{rCl}
\label{eq: Ent08}
\EntCond{X_r}{X_{\set{L}}} \leq \EntCond{X_r}{X_{\Ngb{r}}}.
\end{IEEEeqnarray}
Combining \eqref{eq: Ent07} and \eqref{eq: Ent08} gives that
\begin{IEEEeqnarray}{rCl}
\label{eq: Ent09}
\EntCond{X_{\set{R}_d}}{X_{\set{L}}} \leq \sum_{r \in \set{R}_d} \EntCond{X_r}{X_{\Ngb{r}}}.
\end{IEEEeqnarray}
For $r \in \set{R}_d$, let
\begin{IEEEeqnarray}{rCl}
\label{eq: def RV Q}
Q_r \eqdef \I{ \set{S} \cap \Ngb{r} = \emptyset}
\end{IEEEeqnarray}
be the indicator function of the event where none of the vertices that
are adjacent (in $G$) to the vertex $r$ are included in the (randomly selected)
independent set $\set{S}$. Then,
\begin{IEEEeqnarray}{rCl}
\label{eq: Ent10}
\EntCond{X_r}{X_{\Ngb{r}}} \leq \EntCond{X_r}{Q_r}
\end{IEEEeqnarray}
since the random vector $X_{\Ngb{r}}$ indicates which of the indices $i \in \Ngb{r}$
are included in $\set{S}$, whereas the binary random variable $Q_r$ only indicates
if there is such an index. Consequently, \eqref{eq: Ent09} and \eqref{eq: Ent10}
imply that
\begin{IEEEeqnarray}{rCl}
\label{eq: Ent11}
\EntCond{X_{\set{R}_d}}{X_{\set{L}}} \leq \sum_{r \in \set{R}_d} \EntCond{X_r}{Q_r}.
\end{IEEEeqnarray}
For the binary random variable $Q_r$, let
\begin{IEEEeqnarray}{rCl}
\label{eq: q_r}
q_r \eqdef \prob[Q_r = 1].
\end{IEEEeqnarray}
By \eqref{eq: def RV Q}, $Q_r = 0$ if and only if $\set{S} \cap \Ngb{r} \neq \emptyset$,
which implies that $r \not\in \set{S}$ since there is a vertex in $\Ngb{r}$ that belongs
to the independent set $\set{S}$. Therefore, if $Q_r = 0$, then $X_r = 0$ (see \eqref{eq: indicators}),
so
\begin{IEEEeqnarray}{rCl}
\label{eq: Ent12}
\EntCond{X_r}{Q_r = 0} = 0.
\end{IEEEeqnarray}
If $Q_r = 1$, then $X_r \in \{0,1\}$ and it is also equiprobable (the latter holds since given $Q_r = 1$,
the independent set $\set{S}$ is uniformly distributed over all the independent sets in $\indset{G}$ that
do not include any neighbor of the vertex $r$, so the vertex $r$ can be either removed from or added to
such an independent set, while still giving an independent set that does not include any neighbor of $r$).
Hence,
\begin{IEEEeqnarray}{rCl}
\label{eq: Ent13}
\EntCond{X_r}{Q_r = 1} = 1.
\end{IEEEeqnarray}
Hence, from \eqref{eq: q_r}--\eqref{eq: Ent13},
\begin{IEEEeqnarray}{rCl}
\EntCond{X_r}{Q_r} &=& q_r \, \EntCond{X_r}{Q_r = 1} + (1-q_r) \, \EntCond{X_r}{Q_r = 0} \nonumber \\
&=& q_r,  \label{eq: Ent14}
\end{IEEEeqnarray}
and the combination of \eqref{eq: Ent11} and \eqref{eq: Ent14} yields
\begin{IEEEeqnarray}{rCl}
\label{eq: Ent15}
\EntCond{X_{\set{R}_d}}{X_{\set{L}}} \leq \sum_{r \in \set{R}_d} q_r.
\end{IEEEeqnarray}

We next upper bound $\Ent{X_{\set{L}_d}}$, which is the first summand in the right side
of \eqref{eq: Ent06}, and here Shearer's lemma (see Proposition~\ref{proposition: Shearer's Lemma})
comes into the picture. Since, by definition, $\set{R}_d$ is the set of the vertices that are
connected to the subset $\set{L}_d$ of the degree-$d$ vertices in $\set{L}$, and $\Ngb{r}$ is the set
of vertices in $\set{L}$ that are connected to a vertex $r \in \set{R}_d$ in the bipartite graph $G$,
then it follows that every vertex in $\set{L}_d$ belongs to at least $d$ of the subsets
$\{\Ngb{r}\}_{r \in \set{R}_d}$. Hence, by Shearer's lemma (which, in light of
Remark~\ref{remark: Sheaerer's lemma}, it also holds regardless of the fact
that, for $r \in \set{R}_d$, the set $\Ngb{r}$ is not necessarily a subset of $\set{L}_d$),
\begin{IEEEeqnarray}{rCl}
\label{eq: Ent16}
\Ent{X_{\set{L}_d}} \leq \frac1d \sum_{r \in \set{R}_d} \Ent{X_{\Ngb{r}}}.
\end{IEEEeqnarray}
The binary random variable $Q_r$ is a deterministic function of the random vector
$X_{\Ngb{r}}$ since, from \eqref{eq: indicators} and \eqref{eq: def RV Q}, $Q_r=1$
if and only if all the entries of $X_{\Ngb{r}}$ are equal to~0. Consequently, for
all $r \in \set{R}_d$,
\begin{IEEEeqnarray}{rCl}
\Ent{X_{\Ngb{r}}} &=& \Ent{X_{\Ngb{r}}, Q_r} \label{eq: Ent17} \\[0.1cm]
&=& \Ent{Q_r} + \EntCond{X_{\Ngb{r}}}{Q_r} \label{eq: Ent18} \\
&=& \EntBin{q_r} + \EntCond{X_{\Ngb{r}}}{Q_r},  \label{eq: Ent19}
\end{IEEEeqnarray}
where the equality in \eqref{eq: Ent19} follows from \eqref{eq2: EntBin} and \eqref{eq: q_r}.
Next, from \eqref{eq: q_r},
\begin{IEEEeqnarray}{rCl}
\label{eq: Ent20}
\EntCond{X_{\Ngb{r}}}{Q_r} = q_r \EntCond{X_{\Ngb{r}}}{Q_r=1} + (1-q_r) \EntCond{X_{\Ngb{r}}}{Q_r=0}.
\end{IEEEeqnarray}
If $Q_r = 1$, then $X_{\Ngb{r}}$ is a vector of zeros, so
\begin{IEEEeqnarray}{rCl}
\label{eq: Ent21}
\EntCond{X_{\Ngb{r}}}{Q_r=1} = 0.
\end{IEEEeqnarray}
Otherwise, if $Q_r=0$, then $X_i=1$ for at least one element $i \in \Ngb{r}$; since
$\card{\Ngb{r}}= \dgr{r}$ is the degree of the vertex $r$ (by assumption, there are
no multiple edges connecting any pair of vertices), it follows that the vector
$X_{\Ngb{r}} \in \{0,1\}^{\dgr{r}}$ cannot be the zero vector, so
\begin{IEEEeqnarray}{rCl}
\label{eq: Ent22}
\EntCond{X_{\Ngb{r}}}{Q_r=0} \leq \log(2^{\dgr{r}}-1).
\end{IEEEeqnarray}
Combining \eqref{eq: Ent17}--\eqref{eq: Ent22} gives
\begin{IEEEeqnarray}{rCl}
\label{eq: Ent23}
\Ent{X_{\Ngb{r}}} \leq \EntBin{q_r} + (1-q_r) \log(2^{\dgr{r}}-1),
\end{IEEEeqnarray}
and, from \eqref{eq: Ent16} and \eqref{eq: Ent23}, we get the following upper bound on
the first summand in the right side of \eqref{eq: Ent06}:
\begin{IEEEeqnarray}{rCl}
\label{eq: Ent24}
\Ent{X_{\set{L}_d}} \leq \frac1d \sum_{r \in \set{R}_d}
\bigl\{ \EntBin{q_r} + (1-q_r) \log(2^{\dgr{r}}-1) \bigr\}.
\end{IEEEeqnarray}
Consequently, combining \eqref{eq: Ent02}--\eqref{eq: Ent06}, \eqref{eq: Ent15} and
\eqref{eq: Ent24} implies that
\begin{IEEEeqnarray}{rCl}
\label{eq: Ent25}
\Ent{X^n} &\leq& \sum_{d \in \set{D}_{\mathrm{L}}} \bigl\{ \Ent{X_{\set{L}_d}} + \EntCond{X_{\set{R}_d}}{X_{\set{L}}} \bigr\} \\
\label{eq: Ent26}
&\leq& \sum_{d \in \set{D}_{\mathrm{L}}} \Biggl\{ \, \sum_{r \in \set{R}_d} q_r + \frac1d \sum_{r \in \set{R}_d}
\bigl\{ \EntBin{q_r} + (1-q_r) \log(2^{\dgr{r}}-1) \bigr\} \Biggr\} \\[0.1cm]
\label{eq: Ent27}
&=& \sum_{d \in \set{D}_{\mathrm{L}}} \Biggl\{ \frac1d \sum_{r \in \set{R}_d}
\bigl\{ \EntBin{q_r} + (1-q_r) \log(2^{\dgr{r}}-1) + q_r \log(2^d) \bigr\} \Biggr\} \\[0.1cm]
\label{eq: Ent28}
&=& \sum_{d \in \set{D}_{\mathrm{L}}} \Biggl\{ \frac1d \sum_{r \in \set{R}_d}
\biggl\{ \EntBin{q_r} + q_r \log \biggl(\frac{2^d}{2^{\dgr{r}}-1} \biggr) + \log(2^{\dgr{r}}-1) \biggr\} \Biggr\}.
\end{IEEEeqnarray}
Since $q_r \in [0,1]$ for $r \in \set{R}_d$, we next maximize an auxiliary function
$f_r \colon [0,1] \to \Reals$, defined as
\begin{IEEEeqnarray}{rCl}
\label{eq: f_r}
f_r(x) \eqdef \EntBin{x} + x \log \biggl(\frac{2^d}{2^{\dgr{r}}-1} \biggr), \quad x \in [0,1],
\end{IEEEeqnarray}
in order to obtain an upper bound on the right side of \eqref{eq: Ent28} which is independent of $\{q_r\}$.
By \eqref{eq2: EntBin}, setting the first derivative of the concave function $f_r(\cdot)$ to zero gives the equation
\begin{IEEEeqnarray}{rCl}
\label{derivative f_r is 0}
\log \biggl(\frac{1-x}{x}\biggr) + \log \biggl(\frac{2^d}{2^{\dgr{r}}-1} \biggr) = 0,
\end{IEEEeqnarray}
whose solution is given by
\begin{IEEEeqnarray}{rCl}
\label{eq: solution x}
x = \frac{2^d}{2^d + 2^{\dgr{r}} - 1}.
\end{IEEEeqnarray}
Consequently, it follows from \eqref{eq: Ent25}--\eqref{eq: f_r} and \eqref{eq: solution x} that
\begin{IEEEeqnarray}{rCl}
\label{eq: Ent29}
\Ent{X^n} &\leq& \sum_{d \in \set{D}_{\mathrm{L}}} \Biggl\{ \frac1d \sum_{r \in \set{R}_d}
\Bigl\{ f_r(q_r) + \log(2^{\dgr{r}}-1) \Bigr\} \Biggr\} \\
\label{eq: Ent30}
&\leq& \sum_{d \in \set{D}_{\mathrm{L}}} \Biggl\{ \frac1d \sum_{r \in \set{R}_d}
\Bigl\{ f_r\biggl( \frac{2^d}{2^d + 2^{\dgr{r}} - 1} \biggr) + \log(2^{\dgr{r}}-1) \Bigr\} \Biggr\}
\end{IEEEeqnarray}
and the calculation of the term in the inner sum in the right side of \eqref{eq: Ent30} gives
\begin{IEEEeqnarray}{rCl}
&& \hspace*{-0.5cm} f_r\biggl( \frac{2^d}{2^d + 2^{\dgr{r}} - 1} \biggr) + \log(2^{\dgr{r}}-1) \nonumber \\[0.1cm]
&& = \biggEntBin{\frac{2^d}{2^d + 2^{\dgr{r}} - 1}} + \biggl(\frac{2^d}{2^d + 2^{\dgr{r}} - 1}\biggr) \,
\log \biggl(\frac{2^d}{2^{\dgr{r}}-1} \biggr) + \log(2^{\dgr{r}}-1) \label{eq: Ent31} \\[0.1cm]
&& = - \biggl(\frac{2^d}{2^d + 2^{\dgr{r}} - 1}\biggr) \, \log\biggl( \frac{2^d}{2^d + 2^{\dgr{r}} - 1} \biggr)
-\biggl(\frac{2^{\dgr{r}} - 1}{2^d + 2^{\dgr{r}} - 1}\biggr) \, \log\biggl( \frac{2^{\dgr{r}} - 1}{2^d + 2^{\dgr{r}} - 1} \biggr) \nonumber \\[0.1cm]
&& \hspace*{0.4cm} + \biggl(\frac{2^d}{2^d + 2^{\dgr{r}} - 1} \biggr) \, \log \biggl(\frac{2^d}{2^{\dgr{r}}-1} \biggr) + \log(2^{\dgr{r}}-1) \label{eq: Ent32} \\[0.1cm]
&& = - \biggl(\frac{2^d}{2^d + 2^{\dgr{r}} - 1}\biggr) \, \biggl[ \log\biggl( \frac{2^d}{2^d + 2^{\dgr{r}} - 1} \biggr)
- \log\biggl( \frac{2^d}{2^{\dgr{r}} - 1} \biggr) \biggr] \nonumber \\[0.1cm]
&& \hspace*{0.4cm} -\biggl(\frac{2^{\dgr{r}} - 1}{2^d + 2^{\dgr{r}} - 1}\biggr) \, \log\biggl( \frac{2^{\dgr{r}} - 1}{2^d + 2^{\dgr{r}} - 1} \biggr)
+ \log(2^{\dgr{r}}-1) \label{eq: Ent33} \\[0.1cm]
&& = -\biggl(\frac{2^d}{2^d + 2^{\dgr{r}} - 1}\biggr) \, \log \biggl( \frac{2^{\dgr{r}} - 1}{2^d + 2^{\dgr{r}} - 1} \biggr)
-\biggl(\frac{2^{\dgr{r}} - 1}{2^d + 2^{\dgr{r}} - 1}\biggr) \, \log\biggl( \frac{2^{\dgr{r}} - 1}{2^d + 2^{\dgr{r}} - 1} \biggr) \nonumber \\[0.1cm]
&& \hspace*{0.4cm} + \log(2^{\dgr{r}}-1) \label{eq: Ent34} \\[0.1cm]
&&= -\log \biggl( \frac{2^{\dgr{r}} - 1}{2^d + 2^{\dgr{r}} - 1} \biggr) + \log(2^{\dgr{r}}-1) \label{eq: Ent35} \\[0.1cm]
&&= \log\bigl(2^d + 2^{\dgr{r}} - 1\bigr), \label{eq: Ent36}
\end{IEEEeqnarray}
where \eqref{eq: Ent31} and \eqref{eq: Ent32} hold, respectively, by \eqref{eq: f_r} and \eqref{eq2: EntBin}.
Substituting the equality in \eqref{eq: Ent36} into the upper bound on the entropy in the right side of \eqref{eq: Ent30},
together with \eqref{eq: Ent01}, gives
\begin{IEEEeqnarray}{rCl}
\label{eq: Ent37}
\log \, \bigcard{\indset{G}} &\leq& \sum_{d \in \set{D}_{\mathrm{L}}}
\Biggl\{ \frac1d \sum_{r \in \set{R}_d} \log\bigl(2^d + 2^{\dgr{r}} - 1\bigr) \Biggr\},
\end{IEEEeqnarray}
which, by exponentiation of both sides of \eqref{eq: Ent37}, gives
\begin{IEEEeqnarray}{rCl}
\label{eq1: UB cardinality}
\bigcard{\indset{G}} \leq \prod_{d \in \set{D}_{\mathrm{L}}} \prod_{r \in \set{R}_d} \bigl(2^d + 2^{\dgr{r}} - 1\bigr)^{\frac1d}.
\end{IEEEeqnarray}
The upper bound in the right side of \eqref{eq1: UB cardinality} is in general looser than the bound in Theorem~\ref{theorem: Sah et al., 2019}.
Indeed, to clarify this point, let $\Gamma_{d,d'}$ denote the fraction of vertices in $\set{R}_d$ having degree $d' \in \set{D}_\mathrm{R}$. Then,
\begin{IEEEeqnarray}{rCl}
\prod_{d \in \set{D_{\mathrm{L}}}} \prod_{r \in \set{R}_d}
(2^d + 2^{\dgr{r}} - 1)^{\frac1d}
&=& \prod_{d \in \set{D}_{\mathrm{L}}} \prod_{d' \in \set{D}_{\mathrm{R}}}
(2^d + 2^{d'} - 1)^{\frac{\card{\set{R}_d} \, \Gamma_{d,d'}}{d}} \label{eq2: UB cardinality}  \\
&=& \prod_{d \in \set{D}_{\mathrm{L}}} \prod_{d' \in \set{D}_{\mathrm{R}}}
\Bigl\{ (2^d + 2^{d'} - 1)^{\frac1{d \, d'}} \Bigr\}^{d' \, \card{\set{R}_d} \, \Gamma_{d,d'}}  \label{eq3: UB cardinality}  \\
&\geq& \prod_{(u,v) \in \E{G}} \bigl( 2^{\dgr{u}} + 2^{\dgr{v}} - 1 \bigr)^{\frac1{\dgr{u} \, \dgr{v}}},  \label{eq4: UB cardinality}
\end{IEEEeqnarray}
where \eqref{eq2: UB cardinality} holds since, for all $d \in \set{D}_{\mathrm{L}}$, the number of vertices in
$\set{R}_d$ with degree $d' \in \set{D}_{\mathrm{R}}$ is equal to $\card{\set{R}_d} \, \Gamma_{d,d'}$; finally,
\eqref{eq4: UB cardinality} holds since the number of edges $e = (u,v) \in \E{G}$ that connect left vertices of
degree $d$ and right vertices of degree $d'$ is less than or equal to $d' \, \card{\set{R}_d} \, \Gamma_{d,d'}$
(since $d'$ edges emanate from each such right vertex, but these edges are not necessarily connected to left vertices of degree $d$).
In view of this explanation, there is however an interesting case where the upper bound in the right side of \eqref{eq1: UB cardinality}
and the bound in Theorem~\ref{theorem: Sah et al., 2019} coincide.

Let $G$ be a bipartite graph that is $d$-regular on one side (i.e, one type of its vertices have a fixed degree $d$, and the other type
of vertices can be irregular with arbitrary degrees). Without any loss of generality, one can assume that the left vertices
are $d$-regular (as otherwise, the graph can be flipped without affecting its independent sets, and also the bound in
Theorem~\ref{theorem: Sah et al., 2019} is symmetric in the degrees $\dgr{u}$ and $\dgr{v}$).
In this setting, $\set{L}_d = \set{L}$ and $\set{R}_d = \set{R}$ (recall that, by assumption, there are no isolated vertices).
Consequently, the right side of \eqref{eq1: UB cardinality} is specialized to
\begin{IEEEeqnarray}{rCl}
\label{eq5: UB cardinality}
\bigcard{\indset{G}} \leq \prod_{r \in \set{R}} \bigl(2^d + 2^{\dgr{r}} - 1\bigr)^{\frac1d}.
\end{IEEEeqnarray}
Since there are exactly $\dgr{r}$ edges connecting each vertex $r \in \set{R}$ with vertices in $\set{L}$, and (by the latter assumption) all
of the left vertices in $\set{L}$ are of a {\em fixed degree} $d$, it follows that in this setting, the right side of \eqref{eq5: UB cardinality}
can be rewritten in the form
\begin{IEEEeqnarray}{rCl}
\prod_{r \in \set{R}} \bigl(2^d + 2^{\dgr{r}} - 1\bigr)^{\frac1d} &=& \prod_{r \in \set{R}}
\left( \bigl(2^d + 2^{\dgr{r}} - 1\bigr)^{\frac1{d \dgr{r}}} \right)^{\dgr{r}} \\
\label{eq6: UB cardinality}
&=& \prod_{(u,v) \in \E{G}} \bigl(2^{\dgr{u}} + 2^{\dgr{v}} - 1\bigr)^{\frac1{\dgr{u} \dgr{v}}},
\end{IEEEeqnarray}
which, indeed, shows that the right side of \eqref{eq1: UB cardinality}
and the bound in Theorem~\ref{theorem: Sah et al., 2019} coincide
for bipartite graphs that are regular on one side of the graph (without restricting the other side).

\section{A Variant of the Proof of Zhao's Inequality}
\label{section: An IT proof of Zhao's inequality}

This section suggests a variant of the proof of Zhao's Inequality in \eqref{eq: Zhao's inequality}
(see \cite[Lemma~2.1]{Zhao10}). Although it is somewhat different from the one in \cite{Zhao10},
this forms in essence a reformulation of Zhao’s proof.

Let $G$ be a finite graph, and let $\bigcard{\V{G}} = n$. Label the vertices in the left and
right sides of the bipartite graph $G \times K_2$ (i.e., the bipartite double cover of $G$)
by $\{(i,0)\}_{i=1}^n$ and $\{(i,1)\}_{i=1}^n$, respectively.

Choose independently and uniformly at random two independent sets $\set{S}_0, \set{S}_1 \in \indset{G}$.
For $i \in \OneTo{n}$, let $X_i, Y_i \in \{0,1\}$ be random variables defined as $X_i = 1$
if and only if $i \in \set{S}_0$, and $Y_i = 1$ if and only if $i \in \set{S}_1$. Then, by the
statistical independence and equiprobable selection of the two independent sets from $\indset{G}$,
we have
\begin{IEEEeqnarray}{rCl}
\Ent{X^n, Y^n} &=& \Ent{X^n} + \Ent{Y^n} \label{eq: Ent38} \\
&=& 2 \log \, \bigcard{\indset{G}},   \label{eq: Ent39}
\end{IEEEeqnarray}
where \eqref{eq: Ent38} holds since $X^n = (X_1, \ldots, X_n)$ and $Y^n = (Y_1, \ldots, Y_n)$ are
statistically independent (by construction), and \eqref{eq: Ent39} holds since they both have an
equiprobable distribution over a set whose cardinality is $\bigcard{\indset{G}}$.

Consider the following set of vertices in $G \times K_2$:
\begin{IEEEeqnarray}{rCl}
\label{eq: set S}
\set{S} &\eqdef& \bigl\{ \set{S}_0 \times \{0\} \bigr\} \bigcup \bigl\{ \set{S}_1 \times \{1\} \bigr\} \\[0.1cm]
&=& \bigcup_{i \in \set{S}_0, \, j \in \set{S}_1} \bigl\{ (i,0), (j,1) \bigr\}.
\end{IEEEeqnarray}
The set $\set{S}$ is not necessarily an independent set in $G \times K_2$ since $\bigl( (i,0), (j,1) \bigr) \in \E{G \times K_2}$
for all $i \in \set{S}_0$ and $j \in \set{S}_1$ for which $(i,j) \in \E{G}$ (see \eqref{eq: edges of bipartite double cover}).
We next consider all $(i,j) \in \E{G}$ such that $X_i = Y_j = 1$. To that end, fix an ordering of all the $2^n$ subsets of $\V{G}$,
and let $\set{T} \in \V{G}$ be the first subset in this particular ordering that includes exactly one endpoint of each edge
$(i,j) \in \E{G}$ for which $X_i = Y_j = 1$. Consider the following replacements:
\begin{itemize}
\item If $(i,0) \in \set{S}$ and $i \in \set{T}$, then $(i,0)$ is replaced by $(i,1)$;
\item Likewise, if $(j,1) \in \set{S}$ and $j \in \set{T}$, then $(j,1)$ is replaced by $(j,0)$.
\end{itemize}
Let $\widetilde{\set{S}}$ be the set of new vertices after these possible replacements.
Then, $\widetilde{S} \in \indset{G \times K_2}$ since all adjacent vertices in $\set{S}$ are no longer connected in $\widetilde{\set{S}}$.
Indeed, there is no way that after (say) a vertex $(i,0)$ is replaced by $(i,1)$, there is another replacement
of a vertex $(j,1)$ by $(j,0)$, for some $j$ such that $(i,j) \in \E{G}$; otherwise, that would
mean that $\set{T}$ contains both $i$ and $j$, which is impossible by construction.

Similarly to the way $X^n, Y^n \in \{0,1\}^n$ were defined, let $\widetilde{X}^n, \widetilde{Y}^n \in \{0,1\}^n$ be defined such that,
for all $i \in \OneTo{n}$, $\widetilde{X}_i = 1$ if and only if $(i,0) \in \widetilde{\set{S}}$, and $\widetilde{Y}_i = 1$ if and only if
$(i,1) \in \widetilde{\set{S}}$.

The mapping from $(X^n,Y^n)$ to $(\widetilde{X}^n, \widetilde{Y}^n)$ is injective. Indeed, it is shown to be injective by finding all
indices $(i,j) \in \E{G}$ such that $\widetilde{X}_i = \widetilde{X}_j = 1$ or $\widetilde{Y}_i = \widetilde{Y}_j = 1$, finding the first
subset $\set{T} \in \V{G}$ according to our previous fixed ordering of the $2^n$ subsets of $\V{G}$ that includes exactly one endpoint
of each such edge $(i,j) \in \E{G}$, and performing the reverse operation to return back to $X^n$ and $Y^n$ (e.g., if $(i,j) \in \E{G}$,
$\widetilde{X}_i = \widetilde{X}_j = 1$ and $i \in \set{T}$ while $j \not\in \set{T}$, then $\widetilde{X}_i = 1$ is transformed back
to $Y_i = 1$, and $\widetilde{X}_j = 1$ is transformed back to $X_j = 1$). Consequently, we get
\begin{IEEEeqnarray}{rCl}
\Ent{X^n, Y^n} &=& \Ent{\widetilde{X}^n, \widetilde{Y}^n} \label{eq: Ent40} \\
& \leq & \log \, \bigcard{\indset{G \times K_2}},  \label{eq: Ent41}
\end{IEEEeqnarray}
where \eqref{eq: Ent40} holds by the injectivity of the mapping from $(X^n, Y^n)$ to $(\widetilde{X}^n, \widetilde{Y}^n)$,
and \eqref{eq: Ent41} holds since $\widetilde{S}$ is an independent set in $G \times K_2$, which implies that
$(\widetilde{X}^n, \widetilde{Y}^n)$ can get at most $\bigcard{\indset{G \times K_2}}$ possible values (by definition,
there is a one-to-one correspondence between $\widetilde{\set{S}}$ and $(\widetilde{X}^n, \widetilde{Y}^n)$).
Combining \eqref{eq: Ent38}, \eqref{eq: Ent39}, \eqref{eq: Ent40} and \eqref{eq: Ent41} gives
\begin{IEEEeqnarray}{rCl}
\label{eq: final}
2 \log \, \bigcard{\indset{G}} \leq \log \, \bigcard{\indset{G \times K_2}},
\end{IEEEeqnarray}
which gives \eqref{eq: Zhao's inequality} by exponentiation of both sides of \eqref{eq: final}.

\subsection*{Acknowledgement}
Correspondence with Ashwin Sah and Mehtaab Sawhney (both are currently mathematics graduate
students at MIT), and the constructive comments in the review process are gratefully acknowledged.

\end{document}